\documentclass[12pt]{amsart}

\usepackage {hyperref}

\usepackage{amssymb}

\usepackage{enumerate}

\makeatletter
\@namedef{subjclassname@2010}{%
  \textup{2010} Mathematics Subject Classification}
\makeatother

\newtheorem {theorem} {Theorem}

\newtheorem* {theoremmichael} {Michael Selection Theorem}
\newtheorem* {theoremfanky} {Fan--Kakutani Fixed Point Theorem}

\newtheorem {proposition} [theorem] {Proposition}
\newtheorem {definition} [theorem] {Definition}

\newcommand {\lclass} [2] {\ensuremath{\mathrm L_{#1} \left( #2 \right) }}
\newcommand {\lsclass} [1] {\ensuremath{\mathit l^{#1} }}

\newcommand {\hclass} [2] {\ensuremath{\mathrm H_{#1} \left( #2 \right) }}
\newcommand {\lclassg} [1] {\ensuremath{\mathrm L_{#1}}}

\newcommand {\hclassg} [1] {\ensuremath{\mathrm H_{#1}}}

\newcommand {\BMO} {\ensuremath {\mathrm {BMO}}}

\DeclareMathOperator* {\supp} {supp}

\newcommand {\weightu} {\ensuremath {\mathit u}}
\newcommand {\weightv} {\ensuremath {\mathit v}}
\newcommand {\weightw} {\ensuremath {\mathit w}}

\begin{document}


\baselineskip=17pt


\title [Corona problem with data in ideal spaces of sequences] {Corona problem with data\\in ideal spaces of sequences}

\author [D.~V.~Rutsky] {Dmitry V. Rutsky}
\address {
Steklov Mathematical Institute\\
St. Petersburg Branch\\
Fontanka 27\\
191023 St. Petersburg, Russia}
\email {rutsky@pdmi.ras.ru}

\date{}

\begin{abstract}
Let $E$ be a Banach lattice on $\mathbb Z$ with order continuous norm.
We show that for any function $f = \{f_j\}_{j \in \mathbb Z}$
from the Hardy space $\hclass {\infty} {E}$ such that
$\delta \leqslant \|f (z)\|_E \leqslant 1$ for all $z$ from the unit disk~$\mathbb D$
there exists some solution $g = \{g_j\}_{j \in \mathbb Z} \in \hclass {\infty} {E'}$,
$\|g\|_{\hclass {\infty} {E'}} \leqslant C_\delta$
of the B\'ezout equation $\sum_j f_j g_j = 1$, also known as the vector-valued corona problem
with data in~$\hclass {\infty} {E}$.
\end{abstract}



\maketitle

The (classical) Corona Problem (see, e.~g., \cite [Appendix~3] {nikolsky}) has the following equivalent formulation:
given a finite number of bounded analytic functions
$f = \left\{f_j\right\}_{j = 1}^N \subset \hclassg {\infty}$ on the unit disk~$\mathbb D$,
is the condition $\inf_{z \in \mathbb D} \max_j |f_j (z)| > 0$ sufficient as well as necessary for the existence
of some solutions $g = \left\{g_j\right\}_{j = 1}^N \subset \hclassg {\infty}$
of the B\'ezout equation $\sum_j f_j g_j = 1$?
A positive answer to this problem was established for the first time by L.~Carleson
in~\cite {carleson1962}, and later a relatively simple proof was found by T.~Wolff
(see, e.~g., \cite [Appendix~3] {nikolsky}); we also mention another approach to the proof
based on the theory of analytic multifunctions (see, e.~g., \cite {slodkowski1986}).

These important results set ground for many subsequent developments.
One question to ask is what estimates are possible for the solutions $g$ in terms of the estimates on $f$.
In particular, estimates such as the $\hclass {2} {\lsclass {2}}$ norm of $g$
make it possible to extend these results to infinite sequences $f$.
We formalize this problem somewhat vaguely for now;
see Section~\ref {statementsoftheresults} below for exact definitions.
Notation $\langle f, g\rangle  = \sum_j f_j g_j$ will be used with suitable
sequences $f = \{f_j\}$ and $g = \{g_j\}$.
Let $E$ be a normed space of sequences, and let
$$
E' = \left\{ g\,\,  \mid \,\, |\langle f, g \rangle| < \infty \text { for all $f \in E$}\right\}
$$
be the space of sequences dual to~$E$
with respect to the pairing~$\langle \cdot, \cdot\rangle $.
Suppose that $f \in \hclass {\infty} {E}$ satisfies
$\delta \leqslant \|f (z)\|_{E} \leqslant 1$ for all $z \in \mathbb D$ with some~$\delta > 0$.
We are interested in the existence of a function
$g \in \hclass {\infty} {E'}$ satisfying $\langle f, g\rangle  = 1$, and in the possible norm
estimates of~$g$ in terms of~$\delta$.  If this is the case for all such $f$
then we say that $E$ has the \emph {corona property}.

T.~Wolff's argument allowed M.~Rosenblum, V.~A.~Tolokonnikov and A.~Uchi\-yama to obtain
(independently from one another) a positive
answer to this question in the Hilbert space case, showing that $\lsclass {2}$ has the corona property.
The corresponding estimates were later improved many times; apparently, the best one at the time of this writing can be found
in~\cite {treilwick2005}.
A.~Uchiyama also obtained in~\cite {uchiyama1980} a different
estimate using a
rather involved argument based on the original proof of L.~Carleson, thus establishing that~$\lsclass {\infty}$
also has the corona property.
The intermediate spaces $E = \lsclass {p}$, $2 < p < \infty$ between these two cases can be
reduced to the case $p = 2$
(see~\cite [\S 3] {kisliakovrutsky2012en}),
following Tolokonnikov's (unpublished) remark that Wolff's method can be directly extended to this case at least for
even values of $p > 2$,
but very little was known about the corona property for other spaces~$E$.

Recently, in~\cite {kisliakov2015pen} S.~V.~Kislyakov
extended the corona property
to a large class of Banach ideal sequence spaces~$E$.
Specifically, in this result $E$ is assumed to
be $q$-concave with some $q < \infty$ and satisfy the Fatou property, and
space $\lclass {\infty} {E}$ is assumed to be $\BMO$-regular.
These conditions are satisfied by all $\mathrm {UMD}$ lattices with the Fatou property, and in particular
by spaces~$\lsclass {p}$, $1 \leqslant p < \infty$.
A novel and somewhat counterintuitive idea leading to this result is that for suitable spaces $E_0$ and $E_1$
the corona property of $E_0$ implies the corona property of their pointwise product $E_0 E_1$.
The proof uses the
theory of interpolation for Hardy-type spaces to reduce the result to the well-known case $E = \lsclass {2}$.

In the present work we show how the approach of~\cite {kisliakov2015pen} can be modified to obtain
a complete answer: it turns out that \emph {all} ideal sequence spaces with order continuous norm
have the corona property.  Note that this, in particular, includes all finite-dimensional ideal sequence spaces.
These conditions are more general than the result~\cite {kisliakov2015pen}; see Proposition~\ref {myctmog}
at the end of Section~\ref {cpmpr2}.
It remains unclear if the assumption of the order continuity of the norm can be weakened.

Compared to~\cite {kisliakov2015pen}, our proof relies on somewhat less elementary means, namely
we use a fixed point theorem and a selection theorem to reduce the problem to
Uchiyama's difficult case $E = \lsclass {\infty}$, but otherwise the reduction appears to be
rather simple and straightforward.
Moreover, for $q$-concave lattices $E$ with $q < \infty$ the problem is still reduced in this manner to the relatively
easy standard case $E = \lsclass {2}$ using the same argument as in~\cite {kisliakov2015pen}.

We also mention that very recently in~\cite {zlotnikov2016}
the method described in the present work was also applied to the problem of characterizing the ideals
$I (f) = \left\{ \langle  f, g\rangle  \mid g \in \hclass {\infty} {E'} \right\}$.
Certain classical results
concerning the case $E = \lsclass {2}$ were extended to the case $E = \lsclass {1}$.
The approach~\cite {kisliakov2015pen} based on interpolation of Hardy-type spaces
does not seem to lend itself to such an extension.

We briefly outline the implications for the estimates $C_{E, \delta}$
of the norms of the solutions $g$ in terms of $\delta$.
Theorem~\ref {cpmultt} below can be stated quantitatively: $C_{E_0 E_1, \delta} \leqslant C_{E_0, \frac \delta 2}$ for
$0 < \delta < 1$, provided that $E_0 E_1$ is a Banach lattice with order continuous norm.
Thus for any lattice~$E$ with order continuous norm
we obtain an estimate
\begin {equation}
\label {bascest}
C_{E, \delta} \leqslant C_{\lsclass {\infty}, \frac \delta 2} \leqslant c_1 \delta^{-c_2}
\end {equation}
with
some constants $c_1, c_2 > 0$ independent of $E$.
Furthermore, if a Banach lattice $E$ is $q$-concave with some $1 < q < \infty$ then
$E = \lsclass {q} E_1$ with a Banach lattice $E_1$ (see the proof of Proposition~\ref {myctmog} below), and we get an estimate
$C_{E, \delta} \leqslant C_{\lsclass {q}, \frac \delta 2}$ that may be sharper than~\eqref {bascest} for some values of $q$.
Indeed, we also have an estimate
$C_{\lsclass {p}, \delta} \leqslant C_{\lsclass {2}, \delta^{\frac p 2}}$ for $p \geqslant 2$
(see, e.~g., \cite [\S 3] {kisliakovrutsky2012en})
and $C_{\lsclass {2}, \delta} \leqslant \frac 1 \delta + c \frac 1 {\delta^2} \log \frac 1 \delta$ by~\cite {treilwick2005}
with an explicit constant $c \approx 8.4$.
The latter estimate is known to be close to optimal in terms of the rate of growth as $\delta \to 0$.
Thus, for $q$-concave lattices $E$ with some $2 \leqslant q < \infty$ we also have an estimate
$C_{E, \delta} \leqslant c \frac 1 {\delta^q} (\log \frac 1 \delta)^{\frac q 2}$ for small enough $\delta$.
Our knowledge about sharp estimates for the value of $C_{\lsclass {p}, \delta}$ with $p \neq 2$ seems to be lacking.

\section {Statements of the results}

\label {statementsoftheresults}

A quasi-normed \emph {lattice}
$X$ of measurable functions on a measurable space~$\Omega$, also called an \emph {ideal space},
is a quasi-normed space of measurable functions
such that $f \in X$ and $|g| \leqslant |f|$ implies $g \in X$ and $\|g\|_X \leqslant \|f\|_X$.
Ideal spaces of sequences $E$ are lattices on $\Omega = \mathbb Z$.
A lattice~$X$ is said to have \emph {order continuous quasi-norm} if for any sequence $f_n \in X$
such that $\sup_n |f_n| \in X$ and $f_n \to 0$ almost everywhere one also has $\|f_n\|_X \to 0$.
Lattices~$\lsclass {p}$ have order continuous quasi-norm if and only if $p < \infty$.
For a lattice $X$ of measurable functions the \emph {order dual} $X'$, also called the \emph {associate space},
is the lattice of all measurable functions~$g$
such that the norm
$$
\|g\|_{X'} = \sup_{ f \in X, \, \|f\|_X \leqslant 1}  \int |f g| 
$$
is finite.  For example, the order dual of $\lsclass {p}$ is $\lsclass {p'}$ for all $1 \leqslant p \leqslant \infty$.
If~$X$ is a Banach lattice,
the order dual is contained in the topological dual space~$X^*$ of all continuous linear functionals on~$X$,
and $X' = X^*$ if and only if~$X$ has order continous norm.
For more on lattices see, e.~g.,~\cite {kantorovichold}.

\begin {definition}
\label {cpdef}
Suppose that $E$ is a normed lattice on $\mathbb Z$.
We say that $E$ has the corona property with constant $C_\delta$, $0 < \delta < 1$,
if for any $f \in \hclass {\infty} {E}$ such that\footnote {
Replacing this condition with the nonstrict inequality $\delta < \|f (z)\|_E \leqslant 1$ would allow to
simplify somewhat the arguments in Section~\ref {cpmpred} below; however, the closed form looks nicer.
}
$\delta \leqslant \|f (z)\|_E \leqslant 1$ for all $z \in \mathbb D$
there exists some $g \in \hclass {\infty} {E'}$
such that $\|g\|_{\hclass {\infty} {E'}} \leqslant C_\delta$ and $\langle f (z), g (z)\rangle  = 1$ for all $z \in \mathbb D$.
Such a function~$f$ is called the data for the corona problem with lower bound $\delta$,
and such a function $g$ is called the solution for the corona problem with data $f$.
\end {definition}

For any two quasi-normed lattices $E_0$ and $E_1$ on the same measurable space
the set of pointwise products
$$
E_0 E_1 = \{ h_0 h_1 \mid h_0 \in E_0,\break h_1 \in E_1\}
$$
is a quasi-normed lattice with the quasi-norm defined by
$$
\|h\|_{E_0 E_1} = \inf_{h = h_0 h_1} \|h_0\|_{E_0} \|h_1\|_{E_1}.
$$
For example, the H\"older inequality shows that $\lsclass {p} \lsclass {q} = \lsclass {r}$ with
$\frac 1 r = \frac 1 p + \frac 1 q$.
It is easy to see that $X \lclassg {\infty} = X$ for any lattice $X$.

\begin {theorem}
\label {cpmultt}
Suppose that $E_0$, $E_1$ are Banach lattices on $\mathbb Z$
such that $E = E_0 E_1$ is also a Banach lattice having order continuous norm.
If $E_0$ has the corona property with constants $C_\delta$, $0 < \delta < 1$, then $E$ also
has it with constants $C_{\frac \delta 2}$.
\end {theorem}

The proof of Theorem~\ref {cpmultt} is given in Section~\ref {cpmpr2} below.
Since $\lsclass {\infty}$ has the corona property by~\cite {uchiyama1980}, applying Theorem~\ref {cpmultt} with
$E_0 = \lsclass {\infty}$ and $E_1 = E$ yields the main result, stated as follows.

\begin {theorem}
\label {ecoronat}
Every Banach lattice on $\mathbb Z$ with order continous norm has the corona property.
\end {theorem}

\section {Proof of Theorem~\ref {cpmultt}}

\label {cpmpr2}

We begin with some preparations.
The following result shows that finite-dimensional approximation of the corona property is possible
under the assumption that $E$ has order continous norm.
\begin {proposition}
\label {cpfdred}
Suppose that $E$ is a Banach lattice on $\mathbb Z$ with order continous norm
such that for any $\varepsilon > 0$ and a finite $I \subset \mathbb Z$ the
restriction of $E$ onto $I$ has the corona property with constants
$(1 + \varepsilon) C_\delta$, $0 < \delta < 1$, where $C_\delta$ are independent of~$I$ and~$\varepsilon$.
Then~$E$ has the corona property with constants $C_\delta$, $0 < \delta < 1$.
\end {proposition}
The proof of Proposition~\ref {cpfdred} is given in Section~\ref {cpmpred} below.

A set-valued map $\Phi : X \to 2^Y$ between normed spaces is called \emph {lower semicontinuous}
if for any $x_n \in X$, $x_n \to x$ in $X$ and $y \in \Phi (x)$ there exists a subsequence $n'$ and some
$y_{n'} \in \Phi (x_{n'})$ such that $y_{n'} \to y$ in $Y$.  We need the following well-known result.
\begin {theoremmichael} [{\cite {michael1956}}\footnote {
For the sake of simlicity we omitted the converse part of this famous theorem.
}]
Let $Y$ be a Banach space, $X$ a paracompact space and $\varphi : X \to 2^Y$
a lower semicontinuous multivalued
map taking values that are nonempty, convex and closed.
Then there exists a continuous selection $f : X \to Y$ of $\varphi$, i.~e. $f (x) \in \varphi (x)$ for all $x \in X$.
\end {theoremmichael}

This allows us to conclude that the factorization corresponding to
the product of finite-dimensional Banach lattices can be made continuous
(in the more general infinite-dimensional cases this seems to be unclear).
\begin {proposition}
\label {multcsel}
Suppose that $F_0$ and $F_1$ are finite-dimensional Banach lattices of functions on the same measurable space.
Then for every $\varepsilon$ there exists a continuous map
$\Delta : F_0 F_1 \setminus \{0\} \to F_1$ taking nonnegative values such that
$\|\Delta f\|_{F_1} \leqslant 1$ and
$\|f (\Delta f)^{-1}\|_{F_0} \leqslant (1 + \varepsilon) \|f\|_{F_0 F_1}$.
\end {proposition}
Indeed, we first consider a set-valued map
$\Delta_0 : F_0 F_1 \setminus \{0\} \to 2^{F_1}$
defined~by
\begin {multline*}
\Delta_0 (f) = \left\{ g \in F_1 \mid \text {$g > 0$ everywhere},\right.
\\
\left.
\|g\|_{F_1} < 1,
\left\|f g^{-1}\right\|_{F_0} < (1 + \varepsilon) \|f\|_{F_0 F_1} \right\}
\end {multline*}
for $f \in F_0 F_1 \setminus \{0\}$.
By the definition of the space $F_0 F_1$, the map $\Delta_0$ takes nonempty
values.
It is easy to see that $\Delta_0$ has open graph,
and hence $\Delta_0$ is a lower semicontinuous map.
The graph of a map
$\overline\Delta_0 : F_0 F_1 \setminus \{0\} \to 2^{F_1}$
defined by
\begin {equation*}
\overline\Delta_0 (f) = \left\{ g \in F_1 \mid g \geqslant 0,
\|g\|_{F_1} \leqslant 1,
\left\|f g^{-1}\right\|_{F_0} \leqslant (1 + \varepsilon) \|f\|_{F_0 F_1} \right\}
\end {equation*}
(with the conventions $0 \cdot 0^{-1} = 0$ and $a \cdot 0^{-1} = \infty$ for $a \neq 0$)
is easily seen to be the closure of the graph of the map of $\Delta_0$, therefore
$\overline\Delta_0$ is also a lower semicontinuous map.  The values of $\Delta_0$ are convex and closed,
so by the Michael selection theorem
$\overline\Delta_0$ admits a continuous selection $\Delta$, that is, $\Delta (f) \in \overline\Delta_0 (f)$ for
all $f \in F_0 F_1$.  This selection satisfies the conclusion of Proposition~\ref {multcsel}.

\begin {theoremfanky} [{\cite {fanky1952}}]
Suppose that $K$ is a compact set in a locally convex linear topological space.
Let $\Phi$ be a mapping from $K$ to the set of nonempty subsets of $K$ that are convex and compact, and assume that 
the graph of $\Phi$ is closed.
Then $\Phi$ has a fixed point, i.~e. $x \in \Phi (x)$ for some $x \in K$.
\end {theoremfanky}

A quasi-normed lattice~$X$ of measurable functions is said to have the
\emph {Fatou property} if for any $f_n, f \in X$ such that $\|f_n\|_X
\leqslant 1$ and the sequence~$f_n$ converges to~$f$ almost everywhere
it is also true that $f \in X$ and $\|f\|_X \leqslant 1$.

The following formula (also appearing in \cite [Lemma~1] {kisliakov2015pen}) seems to be rather well known;
see, e.~g., \cite [Theorem~3.7] {schep2010}.
\begin {proposition}
\label {xyxm}
Suppose that $X$ and $Y$ are Banach lattices of measurable functions on the same measurable space
having the Fatou property
such that $X Y$ is also a Banach lattice.  Then $X' = (X Y)' Y$.
\end {proposition}

In order to achieve the best estimate possible with the method used
without assuming that $C_\delta$ is continuous in $\delta$, we take
advantage of the fact that the decomposition in the definition of the
pointwise product of Banach lattices can be made exact if both
lattices satisfy the Fatou property.
\begin {proposition}
\label {fatouexactd}
Let $X$ and $Y$ be Banach lattices of measurable functions on the same measurable space having the Fatou property.
Then for every function $f \in X Y$ there exist some $g \in X$ and $h \in Y$ such that $f = g h$ and
$\|g\|_X \|h\|_Y \leqslant \|f\|_{X Y}$.
\end {proposition}
This appears to be rather well known but hard to find in the literature, so we give a proof.
We may assume that $\|f\|_{X Y} = 1$.
Let $\varepsilon_n \to 0$ be a decreasing sequence.  Sets
\begin {multline*}
F_n = \left\{ g \mid g \geqslant 0, \,\,\text {$\supp g = \supp f$\ up to a set of measure $0$},\right.
\\
\left.
\|g\|_X \leqslant 1, \left\|f g^{-1} \right\|_Y \leqslant 1 + \varepsilon_n \right\} \subset X
\end {multline*}
are nonempty and form a nonincreasing sequence.
It is easy to see that $F_n$ are convex (one uses the convexity of the map $t \mapsto t^{-1}$, $t > 0$).
By the Fatou property of~$X$ and~$Y$ sets~$F_n$ are closed with respect to the convergence in measure,
and they are bounded in~$X$.  The intersection of such a sequence of sets is nonempty
(see~\cite [Chapter 10, \S 5, Theorem~3] {kantorovichold}), so there exists some
$g \in \bigcap_n F_n$, which together with $h = f g^{-1}$ yields the required decomposition.

Now we are ready to prove Theorem~\ref {cpmultt}.
Suppose that under its assumptions $E = E_0 E_1$,
lattice $E_0$ has the corona property with
constant $C_\delta$ for some $0 < \delta \leqslant 1$ and we are given some
$f \in \hclass {\infty} {E}$ such that
$\delta \leqslant \|f (z)\|_E \leqslant 1$ for all $z \in \mathbb D$; we need
to find a suitable $g \in \hclass {\infty} {E'}$ solving $\langle  f, g \rangle  = 1$.

Proposition~\ref {cpfdred} allows us to assume that the lattices have finite support $I \subset \mathbb Z$,
and moreover, we may relax the claimed estimate for the norm of a solution to $(1 + \varepsilon) C_{\frac \delta 2}$
for arbitrary $\varepsilon > 0$.
It is easy to see that finite-dimensional lattices always have the Fatou property, so we may assume that
both $E_0$ and $E_1$, and thus both $\lclass {\infty} {E_0}$ and $\lclass {\infty} {E_1}$
have the Fatou property.  By Proposition~\ref {fatouexactd}
 there exist
some $\weightu \in \lclass {\infty} {E_0}$ and $\weightv \in \lclass {\infty} {E_1}$ such that
$|f| = \weightu \weightv$ and
$\|\weightu\|_{\lclass {\infty} {E_0}} \|\weightv\|_{\lclass {\infty} {E_1}} \leqslant \|f\|_{\lclass {\infty} {E}} \leqslant 1$.
We may further assume that $\|\weightu\|_{\lclass {\infty} {E_0}} \leqslant 1$ and
$\|\weightv\|_{\lclass {\infty} {E_1}} \leqslant 1$.

Since $f = \{f_j\}_{j \in I}$ is analytic and bounded, if we restrict $I$ so that $\mathbb T \times I$ becomes
the support of $f$, we may assume that $\log |f_j| \in \lclassg {1}$ for all $j \in I$.
Boundedness of $\weightu$ and $\weightv$ further implies that $\log |\weightv_j| \in \lclassg {1}$ for all $j \in I$.

Let us fix some $\varepsilon > 0$ and
a sequence $0 < r_j < 1$ such that $r_j \to 1$.  We denote by $P_r$ the operator of convolution with the
Poisson kernel for radius $0 < r < 1$, that is, $P_r a (z) = a (r z)$ for any harmonic function $a$ on $\mathbb D$ and any
$z \in \mathbb D$.

Let
\begin {equation}
\label {bsetdef}
B = \left\{ \log \weightw \mid \weightw \in \lclass {\infty} {E_1}, \|\weightw\|_{\lclass {\infty} {E_1}} \leqslant 2,
\weightw \geqslant \weightv\right\} \subset \lclassg {1}.
\end {equation}
This set is convex, which follows from the well-known logarithmic convexity of the norm of a Banach lattice.
We endow $B$ with the weak topology of~$\lclassg {1}$.
By the Fatou property of $\lclass {\infty} {E_1}$
it is easy to see that $B$ is closed with respect to the convergence in measure,
so $B$ is also closed in $\lclassg {1}$ and thus weakly closed.
The Dunford--Pettis theorem easily shows that $B$ is a compact set, since the functions from $B$
are uniformly bounded from above and below by some summable functions.

For convenience, we denote by $B_Z$ the closed unit ball of a Banach space~$Z$.
We endow $\hclass {\infty} {E_0'}$ with the topology of uniform convergence on compact sets in $\mathbb D \times I$,
and define a (single-valued) map
$
\Phi_0^{(j)} : C_{\frac \delta 2} B_{\hclass {\infty} {E_0'}} \to B
$
by
$$
\Phi_0^{(j)} (h) = \log \left(\Delta (|P_{r_j} h|) + \weightv\right), \quad h \in C_{\frac \delta 2} B_{\hclass {\infty} {E_0'}}
$$
with a map $\Delta$ from Proposition~\ref {multcsel} applied
to $F_0 = E'$ and $F_1 = E_1$ (observe that by Proposition~\ref {xyxm} we have $E_0' = E' E_1$)
and the chosen value of $\varepsilon$.
It is easy to see that $\Phi_0^{(j)}$ is continuous.

We endow 
$\hclass {\infty} {E_1}$ with the topology of uniform convergence on compact sets in $\mathbb D \times I$ and
define a (single-valued) map
$
\Phi_1 : B \to 2 B_{\hclass {\infty} {E_1}}
$
by
\begin {equation}
\label {phi1def}
\Phi_1 (\log \weightw) (z, \omega) =
\exp \left(\frac 1 {2 \pi} \int_0^{2 \pi} \frac {e^{i \theta} + z} {e^{i \theta} - z} \log \weightw \left(e^{i \theta}, \omega\right) d\theta \right)
\end {equation}
for all $\log \weightw \in B$, $z \in \mathbb D$ and $\omega \in I$. 
This map is easily seen to be continuous and
(since the integral under the exponent is the convolution with the Schwarz kernel)
$|\Phi_1 (\log \weightw)| = \weightw$ almost everywhere.

Observe that if $\psi = \Phi_1 (\log \weightw)$ for some $\log \weightw \in B$
and $\varphi = \frac f \psi$ then
$|\varphi| = \frac {|f|} {\weightw} \leqslant \frac {|f|} {\weightv} = \weightu$ and
we have $\varphi \in \hclass {\infty} {E_0}$ with $\|\varphi\|_{\hclass {\infty} {E_0}} \leqslant 1$.
On the other hand,
$$
\delta \leqslant \|f (z)\|_{E} =
\left\|\varphi (z) \psi (z)\right\|_{E_0 E_1} \leqslant \|\varphi (z)\|_{E_0} \|\psi (z)\|_{E_1} \leqslant 2 \|\varphi (z)\|_{E_0},
$$
so $\frac \delta 2 \leqslant \|\varphi (z)\|_{E_0} \leqslant 1$ for all $z \in \mathbb D$.
This means that $\varphi$ belongs to the set
$$
D = \left\{ \varphi \in \hclass {\infty} {E_0} \mid \frac \delta 2 \leqslant \|\varphi (z)\|_{E_0} \leqslant 1 \text { for all $z \in \mathbb D$}\right\}
$$
of corona data functions corresponding to the assumed corona property of~$E_0$.
Thus we may define a (single-valued) map
$$
\Phi_2 : \Phi_1 (B) \to D
$$
by $\Phi_2 (\psi) = \frac f {\psi}$ for $\psi \in \Phi_1 (B)$.
We endow $D$ with the topology of uniform convergence on compact sets in $\mathbb D \times I$.
The continuity of $\Phi_2$ is evident.

We define a set-valued map
$\Phi_3 : D \to 2^{C_{\frac \delta 2} B_{\hclass {\infty} {E_0'}}}$
by
$$
\Phi_3 (\varphi) = \left\{h \in \hclass {\infty} {E_0'} \mid \langle \varphi, h\rangle  = 1,
\|h\|_{\hclass {\infty} {E_0'}} \leqslant C_{\frac \delta 2}
\right\}
$$
for $\varphi \in D$.  By the assumed corona property of $E_0$ map $\Phi_3$ takes nonempty values.
Since the condition $\langle \varphi, h\rangle  = 1$ is equivalent to $\langle \varphi (z), h (z)\rangle  = 1$
for all $z \in \mathbb D$,
it is easy to see that the values of $\Phi_3$ are convex and closed, and thus they are compact.
Similarly, the closedness of the graph of $\Phi_3$ is easily verified.

Now we define a set-valued map $\Phi^{(j)} : C_{\frac \delta 2} B_{\hclass {\infty} {E_0'}} \to 2^{C_{\frac \delta 2} B_{\hclass {\infty} {E_0'}}}$
by
$\Phi^{(j)} = \Phi_3 \circ \Phi_2 \circ \Phi_1 \circ \Phi_0^{(j)}$.  
The graph of $\Phi^{(j)}$ is closed since all individual maps are continuous in the appropriate sense
(specifically, as a composition of upper semicontinous maps, but it is easy to establish the continuity
in this case directly using compactness).
The domain $C_{\frac \delta 2} B_{\hclass {\infty} {E_0'}}$ with the introduced topology
is a compact set in a locally convex linear topological space.  Thus $\Phi^{(j)}$ satisfies the assumptions of the
Fan--Kakutani fixed point theorem, which implies that the maps $\Phi^{(j)}$ admit some fixed points
$h_j \in C_{\frac \delta 2} B_{\hclass {\infty} {E_0'}}$, that is, $h_j \in \Phi^{(j)} (h_j)$ for all $j$.
This means that with $\log \weightw_j = \Phi_0^{(j)} (h_j)$, $\psi_j = \Phi_1 (\log \weightw_j)$ and
$\varphi_j = \Phi_2 (\psi_j)$ we have $h_j \in \Phi_3 (\varphi_j)$.
The first two conditions imply that $|\psi_j| = \Delta (|P_{r_j} h_j|) + v \geqslant \Delta (|P_{r_j} h_j|)$,
so
$$
\left\|(P_{r_j} h_j) (\psi_j)^{-1}\right\|_{\lclass {\infty} {E'}}
\leqslant \left\| \frac {|P_{r_j} h_j|} {\Delta (|P_{r_j} h_j|)}\right\|_{\lclass {\infty} {E'}}
\leqslant \left(1 + \varepsilon\right) C_{\frac \delta 2}
$$
by Proposition~\ref {multcsel}.  Thus
\begin {equation}
\label {phi61}
\left\|\frac {h_j (r_j z)} {\psi_j (z)}\right\|_{E'} \leqslant \left(1 + \varepsilon\right) C_{\frac \delta 2}
\end {equation}
for all $z \in \mathbb D$, and condition $h_j \in \Phi_3 (\varphi_j)$ implies that
\begin {equation}
\label {phi62}
1 = \left\langle \frac {f (z)} {\psi_j (z)}, h_j (z)\right\rangle  = \left\langle f (z), \frac {h_j (z)} {\psi_j (z)}\right\rangle 
\end {equation}
for all $z \in \mathbb D$.
Since sequences $\psi_j$ and $h_j$ are uniformly bounded on compact sets in $\mathbb D \times I$,
by passing to a subsequence we may assume that
$\psi_j \to \psi$ with some $\psi \in \Phi_1 (B)$
and $h_j \to h$ with some $h \in  C_{\frac \delta 2} B_{\hclass {\infty} {E_0'}}$
uniformly on compact sets in $\mathbb D \times I$.  Thus we may pass to the limits
in~\eqref {phi61} and~\eqref {phi62} to see that $\frac h \psi$ is a suitable solution for the corona problem with
data~$f$, which concludes the proof of Theorem~\ref {cpmultt}.

We remark that this construction can be modified to use the Tychonoff fixed point theorem,
which is the particular case of single-valued maps in the setting of the Fan--Kakutani theorem.
It suffices to find a continuous selection for the slightly enlarged map $\Phi_3$, which is the purpose of the next result;
the arbitrarily small increase in the estimate is inconsequential for the scheme of the proof.

\begin {proposition}
\label {cpcsel}
Suppose that a finite-dimensional lattice $E$ has the corona property with constant $C_\delta$ for some $0 < \delta < 1$.
Let
$$
D_E = \left\{ f \in \hclass {\infty} {E} \mid \delta \leqslant \|f (z)\|_{E} \leqslant 1 \text { for all $z \in \mathbb D$}\right\}.
$$
Then for any $\varepsilon > 0$ there exists a continuous map
$$
K : D_E \to (1 + \varepsilon) C_\delta B_{\hclass {\infty} {E'}}
$$
such that $\langle f, K (f)\rangle  = 1$ for any $f \in D_E$.
\end {proposition}
Indeed, let $0 < \alpha < 1$.
We define a set-valued map
$$
K_0 : D_E \to (1 + \alpha) C_\delta B_{\hclass {\infty} {E'}}
$$
by
$$
K_0 (f) = \left\{g \in \hclass {\infty} {E'} \mid \|\langle  f, g \rangle  - 1\|_{\hclassg {\infty}} < \alpha,
\|g\|_{\hclass {\infty} {E'}} < (1 + \alpha) C_\delta
\right\}
$$
for $f \in D_E$.  By the corona property assumption on $E$ map $K_0$ takes nonempty values.
It is easy to see that $K_0$ has open graph and thus $K_0$ is a lower semicontinous map.
The graph of a map
$\overline K_0 : D_E \to (1 + \alpha) C_\delta B_{\hclass {\infty} {E'}}$
defined by
\begin {multline*}
\overline K_0 (f) = \left\{g \in \hclass {\infty} {E'} \mid \|\langle  f, g \rangle  - 1\|_{\hclassg {\infty}} \leqslant \alpha,
\|g\|_{\hclass {\infty} {E'}} \leqslant (1 + \alpha) C_\delta
\right\}
\end {multline*}
is easily seen to be the closure of the graph of the map of $K_0$, and hence
$\overline K_0$ is also a lower semicontinuous map.  The values of $K_0$ are convex and closed.
By the Michael selection theorem there exists a continuous selection~$K_1$
of the map $\overline K_0$, that is, $K_1 (f) \in \overline K_0 (f)$ for
all $f \in D_E$.
Now observe that $|\langle f (z), K_1 (f) (z)\rangle  - 1| \leqslant \alpha$ implies
$|\langle f (z), K_1 (f) (z)\rangle | \geqslant 1 - \alpha$ for all $z \in \mathbb D$ and $f \in D_E$,
so we may set $K (f) = \frac {K_1 (f)} {\langle f, K_1 (f)\rangle }$
and have $\langle f, K (f)\rangle  = 1$ with
$\|K (f)\|_{\hclass {\infty} {E'}} \leqslant \frac {1 + \alpha} {1 - \alpha} C_\delta$.
Choosing $\alpha$ small enough yields the claimed range of~$K$.

Finally, we mention that Theorem~\ref {ecoronat} includes the result \cite [Corollary~2] {kisliakov2015pen}.
This is implied by the following known observation; we give a proof for convenience.
\begin {proposition}
\label {myctmog}
Suppose that~$X$ is a Banach lattice of measurable functions having the Fatou property and $X$ is $q$-concave with some
$1 < q < \infty$.
Then~$X$ has order continuous norm.
\end {proposition}
Lattice $Z^\delta$ is defined by the norm $\|f\|_{Z^\delta} = \left\| |f|^{\frac 1 \delta} \right\|_Z^\delta$
for a quasi-normed lattice $Z$ of measurable functions and $\delta > 0$.
Lattice~$X'$ is $q'$-convex, and hence $Y = (X')^{q'}$ is a Banach lattice with the Fatou property.
Then $X' = Y^{\frac 1 {q'}}$.  The Fatou property is equivalent to the order reflexivity $X = X''$,
and using the well-known formula for the duals of the Calder\'on-Lozanovsky products (see, e.~g., \cite [Theorem~2.10] {schep2010})
we may write
$$
X = (X')' = \left({Y}^{\frac 1 {q'}} {\lclassg {\infty}\strut}^{\frac 1 q}\right)' = Y'^{\frac 1 {q'}} {\lclassg {1}\strut}^{\frac 1 q} = Y'^{\frac 1 {q'}} \lclassg {q}.
$$
Since lattice $\lclassg {q}$ has order continuous norm, it suffices to establish the following.
\begin {proposition}
\label {ocnmult}
Suppose that $X$ and $Y$ are quasi-normed lattices of measurable functions and $Y$ has order continous quasi-norm.
Then $X Y$ also has order continuous quasi-norm.
\end {proposition}
Let $f_n \in X Y$ be a sequence with $f = \sup_n |f_n| \in X Y$ such that $f_n \to 0$ almost everywhere.
Then $f = g h$ with some $g \in X$ and $h \in Y$.
We may assume that $g, h \geqslant 0$.  Sequence $h_n = \frac {f_n} g$ also converges to $0$ almost everywhere,
and $|h_n| \leqslant \frac f g = h$, so $\sup_n |h_n| \in Y$.  By the order continuity of the quasi-norm of $Y$
we have $\|h_n\|_Y \to 0$, hence $\|f_n\|_{X Y} \leqslant \|g\|_X \|h_n\|_Y \to 0$.

\section {Proof of Proposition~\ref {cpfdred}}

\label {cpmpred}

First, observe that if a lattice $E$ on $\mathbb Z$ has order continuous norm
then $\lclass {1} {E}$ also has order continous norm, we have
$\lclass {\infty} {E'} = \left[\lclass {1} {E}\right]' = \left[\lclass {1} {E}\right]^*$,
and $\hclass {\infty} {E'}$ is easily seen to be $w^*$-closed in $\lclass {\infty} {E'}$
(see, e.~g., \cite [\S 1.2.1] {kisliakovrutsky2012en}).
The $w^*$-convergence of a sequence $h_k \in \hclass {\infty} {E'}$ to some $h$ implies that
$h_k (z) \to h (z)$ in the $*$-weak topology of $E' = E^*$ for all $z \in \mathbb D$.

Now suppose that under the assumptions of Proposition~\ref {cpfdred}
$0 < \delta < 1$ and $f \in \hclass {\infty} {E}$ satisfies
$\delta \leqslant \|f (z)\|_E \leqslant 1$ for all $z \in \mathbb D$.
Let $I_k \subset \mathbb Z$ be a nondecreasing sequence such that $\bigcup_k I_k = \mathbb Z$,
and fix a sequence $\varepsilon_j > 0$, $\varepsilon_j \to 0$.
We consider the natural approximations
$f_{A, r, k} (z) = A f (r z) \chi_{I_k}$, $z \in \mathbb D$, for $0 < r \leqslant 1$ and $A \geqslant 1$.
If there exists some sequence of parameters $A_j \to 1$, $r_j \to 1$ and $k_j \to \infty$ such that
$f_j = f_{A_j, r_j, k_j}$ is a data for the corona problem with lower bound $\delta$ then
by the assumptions there exist
some $g_j$ such that $\langle f_j, g_j\rangle  = 1$ and
$\|g_j\|_{\hclass {\infty} {E'}} \leqslant (1 + \varepsilon_j) C_\delta$.
By passing to a subsequence we may assume that $A_j$ is nonincreasing, $r_j$ and $k_j$ are nondecreasing, and
$g_j \to g$ in the $*$-weak topology of
$\lclass {\infty} {E'}$ for some $g \in \hclass {\infty} {E'}$.
Observe that
\begin {multline}
\label {cpfdreddec}
1 = \langle f_j (z), g_j (z)\rangle  =
\\
\langle f (z), g_j (z)\rangle  +
\langle f (r_j z) - f (z), g_j (z)\rangle  + \langle f_j (z) - f (r_j z), g_j (z)\rangle 
\end {multline}
for all $z \in \mathbb D$.
The first term in~\eqref {cpfdreddec} converges to $\langle f (z), g (z)\rangle $.
Since the $E$-valued analytic function $f$ is strongly continuous at every $z \in \mathbb D$,
the second term in~\eqref {cpfdreddec} converges to $0$.
By the assumptions
\begin {multline}
\label {cpfdreddec2}
|f_j (z) - f (r_j z)| = |A_j f (r_j z) \chi_{I_{k_j}} - f (r_j z)|
\leqslant
\\
\chi_{I_{k_j}} |A_j f (r_j z) \chi_{I_{k_j}} - f (r_j z)| + 
\chi_{\mathbb Z \setminus I_{k_j}} |A_j f (r_j z) \chi_{I_{k_j}} - f (r_j z)|
=
\\
\chi_{I_{k_j}} |A_j f (r_j z) - f (r_j z)| +
\chi_{\mathbb Z \setminus I_{k_j}} |f (r_j z)|
\leqslant
\\
(A_j - 1) |f (r_j z)| +
\chi_{\mathbb Z \setminus I_{k_j}} |f (z)| + |f (r_j z) - f (z)|.
\end {multline}
The norm in $E$ of the first term in~\eqref {cpfdreddec2} is estimated by $(A_j - 1)$, and thus it converges to $0$.
The second term in~\eqref {cpfdreddec2} converges to $0$ in $E$ by the assumption that $E$ has order continous norm.
The third term in~\eqref {cpfdreddec2} converges to $0$ in $E$ by the strong continuity of $f$ in $\mathbb D$.
It follows that the third term in~\eqref {cpfdreddec} is dominated by
$\|f_j (z) - f (r_j z)\|_E \|g_j (z)\|_{E'} \leqslant \|f_j (z) - f (r_j z)\|_E (1 + \varepsilon_j) C_\delta$,
and so it also converges to $0$.
Therefore passing to the limit in \eqref {cpfdreddec} yields $\langle f (z), g (z)\rangle  = 1$ for all $z \in \mathbb D$.
We also have $\|g\|_{\hclass {\infty} {E'}} \leqslant \limsup_j \|g_j\|_{\hclass {\infty} {E'}} \leqslant C_\delta$,
so~$g$ is a solution for the corona problem with data~$f$ having the claimed constant~$C_\delta$.

Thus it suffices to find a suitable sequence of parameters.  We consider two cases.
In the first case
$\left\|f (z)\right\|_E = 1$ for all $z \in \mathbb D$.
We take $A_j = 1$ and any increasing sequence $r_j \to 1$.
By the order continuity of norm we have
$\left\|f (z) \chi_{I_k}\right\|_E \to \left\|f (z)\right\|_E = 1$ for every $z \in \mathbb D$,
so by the compactness of closed sets in $\mathbb D$ and the assumption that $\delta < 1$
we have $\left\|f (r_j z) \chi_{I_k}\right\|_E \geqslant \delta$
for large enough $k_j$.
Thus $f_j$ is a suitable corona data in this case.

In the second case $\left\|f (z_0)\right\|_E < 1$ for some $z_0 \in \mathbb D$.
With the help of an automorphism we may assume for convenience that $z_0 = 0$.
We also fix any increasing sequence $r_j \to 1$.
A simple consequence of the Schwarz lemma (see, e.~g., \cite [Chapter~1, Corollary~1.3] {garnett1981baf}) shows that
\begin {equation}
\label {cpfdreddec3}
|\langle f (z), e'\rangle | \leqslant \frac {|\langle f (0), e'\rangle | + |z|} {1 + |\langle f (0), e'\rangle |\,\, |z|}
\end {equation}
for all $z \in \mathbb D$ and $e' \in E' = E^*$ with $\|e'\|_{E'} \leqslant 1$.
Since function $(x, y) \mapsto \frac {x + y} {1 + x y}$ is increasing in both $x \in [0, 1]$ and $y \in [0, 1]$,
taking the supremum in \eqref {cpfdreddec3} over all such $e'$ yields
$\|f (z)\|_E \leqslant \frac {\|f (0)\|_E + |z|} {1 + \|f (0)\|_E |z|}$,
thus $\alpha_j = \sup_{z \in \mathbb D} \|f (r_j z)\|_E < 1$.
Setting $A_j = \frac 1 {\alpha_j}$ yields $\|A_j f (r_j z)\|_E \leqslant 1$ for all $z \in \mathbb D$.
Again, since $E$ has order continuous norm we have
$\left\|A_j f (z) \chi_{I_k}\right\|_E \to \left\|A_j f (z)\right\|_E \geqslant A_j \delta > \delta$
for every $z \in \mathbb D$, and we also have
$\left\|A_j f (r_j z) \chi_{I_{k_j}}\right\|_E \geqslant \delta$
for large enough $k_j$, so $f_j$ is a suitable corona data in this case as well.
The proof of Proposition~\ref {cpfdred} is complete.

\subsection* {Acknowledgements}

The author is grateful to S.~V.~Kislyakov for stimulating discussions and a surprising conjecture
that eventually became the statement of Theorem~\ref {cpmultt}, and to the referee for thorough and helpful remarks.

\normalsize
\baselineskip=17pt

\bibliographystyle {acmx}

\bibliography {bmora}

\begin{thebibliography}{10}

\bibitem{carleson1962}
{\sc L.~Carleson}.
\newblock Interpolations by bounded analytic functions and the corona problem.
\newblock {\em Ann. Math. 76}, 2 (1962), 547--559.

\bibitem{fanky1952}
{\sc K.~Fan}.
\newblock Fixed-point and minimax theorems in locally convex topological linear
  spaces.
\newblock {\em Proc. Nat. Acad. Sci. U.S.A. 38\/} (1952), 121--126.

\bibitem{garnett1981baf}
{\sc J.~B. Garnett}.
\newblock {\em {Bounded analytic functions}}.
\newblock Academic Press New York, 1981.

\bibitem{kantorovichold}
{\sc L.~V. Kantorovich, G.~P. Akilov}.
\newblock {\em Functional Analysis}, 2nd~ed.
\newblock Nauka, Moscow, 1977.
\newblock In Russian; English transl.: Pergamon Press, Oxford--New York, 1982.

\bibitem{kisliakov2015pen}
{\sc S.~V. Kislyakov}.
\newblock {Corona theorem and interpolation}.
\newblock {\em {Algebra i Analiz} 27}, 5 (2015), 69--80.
\newblock In Russian; English transl.: St. Petersburg Math. J., 27 (2016),
  757--764.

\bibitem{kisliakovrutsky2012en}
{\sc S.~V. Kislyakov, D.~V. Rutsky}.
\newblock Some remarks to the corona theorem.
\newblock {\em Algebra i Analiz 24}, 2 (2012), 171--191.
\newblock In Russian; English transl.: St. Petersburg Math. J., 24 (2013),
  313--326.

\bibitem{michael1956}
{\sc E.~Michael}.
\newblock {Continuous selections. I}.
\newblock {\em Ann. Math. Second Series 63}, 2 (1956), 361--382.

\bibitem{nikolsky}
{\sc N.~K. Nikol'ski\u{i}}.
\newblock {\em Treatise on the Shift Operator}.
\newblock Springer-Verlag, Berlin, Heidelberg, New York, 1986.

\bibitem{schep2010}
{\sc A.~R. Schep}.
\newblock {Products and factors of Banach function spaces}.
\newblock {\em Positivity 14\/} (2010), 301--319.

\bibitem{treilwick2005}
{\sc S.~Treil, B.~D. Wick}.
\newblock {The matrix-valued $H^p$ corona problem in the disk and polydisk}.
\newblock {\em J. Funct. Anal. 226}, 1 (2005), 138--172.

\bibitem{uchiyama1980}
{\sc A.~Uchiyama}.
\newblock {Corona theorems for countably many functions and estimates for their
  solutions}.
\newblock {Preprint, UCLA}, 1980.

\bibitem{slodkowski1986}
{\sc {Z. S\l{}odkowski}}.
\newblock {An analytic set-valued selection and its applications to the corona
  theorem, to polynomial hulls and joint spectra}.
\newblock {\em Trans. Am. Math. Soc. 294}, 1 (1986), 367--377.

\bibitem{zlotnikov2016}
{\sc I.~K. Zlotnikov}.
\newblock {Estimates in problem of ideals in the algebra $H^\infty$}.
\newblock {\em Zap. Nauchn. Sem. POMI 447\/} (2016), 66--74.
\newblock In Russian.

\end{thebibliography}

\end{document}